\documentclass[12pt]{amsart}
\usepackage{amsmath,amssymb,amsbsy,amsfonts,latexsym,amsopn,amstext,
                                               amsxtra,euscript,amscd,bm}
                    
\usepackage{url}
\usepackage[colorlinks,linkcolor=blue,anchorcolor=blue,citecolor=blue]{hyperref}
\usepackage{color}

\begin{document}
\newtheorem{problem}{Problem}
\newtheorem{theorem}{Theorem}
\newtheorem{lemma}[theorem]{Lemma}
\newtheorem{claim}[theorem]{Claim}
\newtheorem{cor}[theorem]{Corollary}
\newtheorem{prop}[theorem]{Proposition}
\newtheorem{definition}{Definition}
\newtheorem{question}[theorem]{Question}

\def\cA{{\mathcal A}}
\def\cB{{\mathcal B}}
\def\cC{{\mathcal C}}
\def\cD{{\mathcal D}}
\def\cE{{\mathcal E}}
\def\cF{{\mathcal F}}
\def\cG{{\mathcal G}}
\def\cH{{\mathcal H}}
\def\cI{{\mathcal I}}
\def\cJ{{\mathcal J}}
\def\cK{{\mathcal K}}
\def\cL{{\mathcal L}}
\def\cM{{\mathcal M}}
\def\cN{{\mathcal N}}
\def\cO{{\mathcal O}}
\def\cP{{\mathcal P}}
\def\cQ{{\mathcal Q}}
\def\cR{{\mathcal R}}
\def\cS{{\mathcal S}}
\def\cT{{\mathcal T}}
\def\cU{{\mathcal U}}
\def\cV{{\mathcal V}}
\def\cW{{\mathcal W}}
\def\cX{{\mathcal X}}
\def\cY{{\mathcal Y}}
\def\cZ{{\mathcal Z}}

\def\A{{\mathbb A}}
\def\B{{\mathbb B}}
\def\C{{\mathbb C}}
\def\D{{\mathbb D}}
\def\E{{\mathbb E}}
\def\F{{\mathbb F}}
\def\G{{\mathbb G}}
\def\I{{\mathbb I}}
\def\J{{\mathbb J}}
\def\K{{\mathbb K}}
\def\L{{\mathbb L}}
\def\M{{\mathbb M}}
\def\N{{\mathbb N}}
\def\O{{\mathbb O}}
\def\P{{\mathbb P}}
\def\Q{{\mathbb Q}}
\def\R{{\mathbb R}}
\def\S{{\mathbb S}}
\def\T{{\mathbb T}}
\def\U{{\mathbb U}}
\def\V{{\mathbb V}}
\def\W{{\mathbb W}}
\def\X{{\mathbb X}}
\def\Y{{\mathbb Y}}
\def\Z{{\mathbb Z}}

\def\ep{{\mathbf{e}}_p}
\def\em{{\mathbf{e}}_m}
\def\eq{{\mathbf{e}}_q}

\def\scr{\scriptstyle}
\def\\{\cr}
\def\({\left(}
\def\){\right)}
\def\[{\left[}
\def\]{\right]}
\def\<{\langle}
\def\>{\rangle}
\def\fl#1{\left\lfloor#1\right\rfloor}
\def\rf#1{\left\lceil#1\right\rceil}
\def\le{\leqslant}
\def\ge{\geqslant}
\def\eps{\varepsilon}
\def\mand{\qquad\mbox{and}\qquad}
\def\tz{\widetilde z}
\def\tX{\widetilde X}
\def\tx{\widetilde x}
\def\tu{\widetilde u}
\def\tv{\widetilde v}
\def\tw{\widetilde w}
\def\tY{\widetilde Y}

\def\sssum{\mathop{\sum\ \sum\ \sum}}
\def\ssum{\mathop{\sum\, \sum}}
\def\ssumw{\mathop{\sum\qquad \sum}}

\def\vec#1{\mathbf{#1}}
\def\inv#1{\overline{#1}}
\def\num#1{\mathrm{num}(#1)}
\def\dist{\mathrm{dist}}

\def\fA{{\mathfrak A}}
\def\fB{{\mathfrak B}}
\def\fC{{\mathfrak C}}
\def\fU{{\mathfrak U}}
\def\fV{{\mathfrak V}}

\newcommand{\bflambda}{{\boldsymbol{\lambda}}}
\newcommand{\bfxi}{{\boldsymbol{\xi}}}
\newcommand{\bfrho}{{\boldsymbol{\rho}}}
\newcommand{\bfnu}{{\boldsymbol{\nu}}}
\newcommand{\Mod}[1]{\ (\mathrm{mod}\ #1)}

\def\GL{\mathrm{GL}}
\def\SL{\mathrm{SL}}

\def\Hba{\overline{\cH}_{a,m}}
\def\Hta{\widetilde{\cH}_{a,m}}
\def\Hb1{\overline{\cH}_{m}}
\def\Ht1{\widetilde{\cH}_{m}}

\def\flp#1{{\left\langle#1\right\rangle}_p}
\def\flm#1{{\left\langle#1\right\rangle}_m}
\def\dmod#1#2{\left\|#1\right\|_{#2}}
\def\dmodq#1{\left\|#1\right\|_q}

\def\Zm{\Z/m\Z}

\def\Err{{\mathbf{E}}}

\newcommand{\comm}[1]{\marginpar{%
\vskip-\baselineskip 
\raggedright\footnotesize
\itshape\hrule\smallskip#1\par\smallskip\hrule}}

\def\xxx{\vskip5pt\hrule\vskip5pt}


\title[Points on polynomial curves in small boxes]{Points on polynomial curves in small boxes modulo an integer}

 \author[B. Kerr] {Bryce Kerr}

\address{Department of Pure Mathematics, University of New South Wales,
Sydney, NSW 2052, Australia}
\email{bryce.kerr89@gmail.com}

 \author[A. Mohammadi] {Ali Mohammadi}

\address{School of Mathematics and Statistics, University of
Sydney, NSW 2006, Australia}
\email{alim@maths.usyd.edu.au}

\date{\today}
\pagenumbering{arabic}


\begin{abstract} 
Given an integer $q$ and a polynomial $f\in \Z_{q}[X]$ of degree $d$ with coefficients in the residue ring $\Z_q=\Z/q\Z,$ we obtain new results concerning the number of solutions to congruences of the form
$$y\equiv f(x) \pmod{q},$$
with integer variables lying in some cube $\cB$ of side length $H$. Our argument uses ideas of Cilleruelo, Garaev, Ostafe and Shparlinski which reduces the problem to the Vinogradov mean value theorem and a lattice point counting problem. We treat the lattice point problem differently using transference principles from the Geometry of Numbers. We also use a variant of the main conjecture for the Vinogradov mean value theorem of Bourgain, Demeter and Guth and of Wooley which allows one to deal with rather sparse sets.
 \end{abstract}
\maketitle
\section{Introduction}
Given an integer $q$ and a polynomial $f\in \Z_{q}[X]$ of degree $d$ with coefficients in the residue ring $\Z_q=\Z/q\Z,$ we consider the problem of estimating the number of solutions to congruences of the form
\begin{align}
\label{eq:mainequation}
y\equiv f(x) \pmod{q},
\end{align}
with integer variables
\begin{align*}
(x,y)\in (K,K+H]\times (L,L+H],
\end{align*}
lying in some cube of side length $H$. This problem and its variants have been considered by a number of previous authors, see for example~\cite{ChSh,Chang,CCGHS,CiGa,CSh,CGOS,Ke,Zh,Zu}. A common strategy for dealing with such equations is to lift~\eqref{eq:mainequation} to a polynomial equation over $\Z$ to which one may apply the following result of Bombieri and Pila~\cite{BomPi}.
\begin{theorem}
\label{thm:bombieripila}
Let $\cC$ be an absolutely irreducible curve of degree $d\ge2$ and suppose $H\ge e^{d^{6}}$. The number of integral points on $\cC$ and inside a square $[1,H]\times [1,H]$  is bounded by
$$O\left(H^{1/d}\exp{(12(\log{H}\log\log{H})^{1/2})}\right).$$
\end{theorem}
To see how one may lift the equation~\eqref{eq:mainequation} to an equation over $\Z$, we first note that by shifting variables and modifying the coefficients of $f$ we may assume $K=L=0$. Now, suppose $f$ is given by
$$f(x)=a_0+a_1x+\dots+a_dx^{d},$$
with leading coefficient satisfying $(a_d,q)=1$. We define the lattice 
\begin{equation*}
\begin{split}
\cL= \{&(z,w_1,\dots,w_d)\in \Z^{d+1} :\\ &\exists \ n\in \Z \ \text{such that} \ n\equiv z \pmod{q}, \ \ a_in \equiv w_i \pmod{q} \},
\end{split}
\end{equation*}
and the convex body
$$D= \left \{ (z,w_1,\dots,w_d)\in \Z^{d+1} \ : \ |z|\le \frac{q}{H}, \ \ |w_i|\le \frac{q}{H^{i}} \right \}.$$
By Minkowski's first theorem,  $D$ contains a non-zero lattice point of $\cL$ provided
\begin{equation}
\label{eq:mft}
\text{Vol}(D)\ge 2^{d+1}\text{Vol}\left(\R^{d+1}/\cL \right).
\end{equation}
Since 
$$\text{Vol}(\R^{d+1}/\cL)=q^{d},$$
and 
$$\text{Vol}(D)=\frac{2^{d+1}q^{d+1}}{H^{(d^2+d+2)/2}},$$
we see that~\eqref{eq:mft} is satisfied provided
\begin{align}
\label{eq:Hcond}
H\le q^{2/(d^2+d+2)}.
\end{align}
Suppose $H$ satisfies~\eqref{eq:Hcond}, let $(z,w_1,\dots,w_d)$ be a nonzero lattice point of $\cL\cap D$ and suppose $n$ is defined by
$$n\equiv z \pmod{q}, \quad a_in\equiv w_i \pmod{q}.$$
If $w_0$ denotes the smallest residue of $a_0n \pmod{q}$ then the congruence
$$y\equiv f(x) \pmod{q},$$
is equivalent to
\begin{align*}
zy \equiv w_0+w_1x+\dots+w_d x^{d} \pmod{q},
\end{align*}
and hence each point satisfying~\eqref{eq:mainequation} must also satisfy
\begin{align*}
zy=w_0+w_1x+\dots+w_d x^{d}+tq,
\end{align*}
for some $t\in \Z$. Since $(a_d,q)=1$ we have $w_d\neq 0$. From $(z,w_1,\dots,w_d)\in D$ we see that there are $O(1)$  possible values for $t$ when the variables $x,y$ satisfy
$$1\le x \le H \quad \text{and} \quad 1\le y \le H,$$
 and for each such value of $t$ we apply Theorem~\ref{thm:bombieripila}. This results in the following.
\begin{theorem}
\label{thm:prelim}
Let $q$ be an integer, $f\in \Z_{q}[X]$ a polynomial of degree $d\ge 2$ with leading coefficient  coprime to $q$ and suppose $\cB$ is a cube of side length $H$. If 
$$H\le q^{2/(d^2+d+2)},$$
then 
$$|\{ (x,y)\in \cB \ : \ y\equiv f(x) \pmod{q}\}|\le H^{1/d+o(1)}.$$
\end{theorem}
In this paper we consider the problem of determining the largest number $\alpha$ depending only on $d$ such that for all integers $q$ and polynomials $f\in \Z_{q}[X]$ of degree $d$ with leading coefficient coprime to $q$ we have the following estimate
\begin{align}
\label{eq:goal}
|\{ (x,y)\in \cB \ : \ y\equiv f(x) \pmod{q}\}|\le H^{1/d+o(1)},
\end{align}
uniformly over all cubes $\cB$ of side length $H\le q^{\alpha}$. When $q$ is prime, improvements on Theorem~\ref{thm:prelim} are known. For example,  Cilleruelo, Garaev, Ostafe and Shparlinski~\cite[Theorem~3]{CGOS} have shown that~\eqref{eq:goal} holds provided $$H\le q^{2/(d^2+3)}.$$ 

 We obtain  new results concerning estimates for the number of solutions to equations of the form~\eqref{eq:mainequation}. In particular, we improve on the above result of Cilleruelo, Garaev, Ostafe and Shparlinski~\cite[Theorem~3]{CGOS} and some results of Chang,  Cilleruelo, Garaev,  Hernández, Shparlinski and  Zumalac\'{a}rregui~\cite[Corollary~3]{CCGHS}. Our argument uses ideas from ~\cite{CCGHS} and~\cite{CGOS}  which reduces the problem of bounding the number of solutions to~\eqref{eq:mainequation} to the Vinogradov mean value theorem and a problem concerning the number of lattice points in the intersection of a convex body. We treat the lattice point problem differently using transference principles from the Geometry of Numbers. We also use some recent progress concerning the main conjecture for the Vinogradov mean value theorem first obtained by Bourgain, Demeter and Guth~\cite{BDG} and shortly after by Wooley~\cite{Wo} which allows one to deal with solutions to the Vinogradov mean value theorem when the variables belong to rather sparse sets, see Lemma~\ref{lem:vmvt} below.

 Finally, we mention techniques from the geometry of numbers have found applications to a wide range of other problems concerning equations in finite fields and residue rings to which we refer the reader to~\cite{BC,BG,BG1,BGKS,Kon}.
\section{Main results}
Our first result improves on a bound of Cilleruelo, Garaev, Ostafe and Shparlinski~\cite[Theorem~3]{CGOS}.
\begin{theorem}
\label{thm:main1}
Let $q$ be an arbitrary integer, $f\in \Z_{q}[X]$ a polynomial of degree $d\ge 2$ with leading coefficient $a_d$ satisfying $(a_d,q)=1$ and suppose $\cB$ is a cube of side length $H$. We have
$$|\{ (x,y)\in \cB \ : \ y\equiv f(x) \pmod{q}\}|\le \frac{H^{1+2/d(d+1)+o(1)}}{q^{2/d(d+1)}}+H^{1/d+o(1)}.$$
In particular, if 
$$H\le q^{2/(d^2+1)},$$
then  
\begin{align}
\label{eq:diagonal}
|\{ (x,y)\in \cB \ : \ y\equiv f(x) \pmod{q}\}|\le H^{1/d+o(1)}.
\end{align}
\end{theorem}
For comparison with Theorem~\ref{thm:main1} we note that~\cite[Theorem~3]{CGOS} gives the bound~\eqref{eq:diagonal} in the shorter range $H\le q^{2/(d^2+3)}.$

Our second result improves on some results of~\cite{CCGHS}.
\begin{theorem}
\label{thm:ConcentEllipticWiderRange}
Let $q$ be an arbitrary integer, $f\in \Z_q[X]$ of degree $d =3$ and leading coefficient $a_3$ satisfying $(a_3,q)=1$ and $\cB$ a cube of side length $H$. Then we have
$$|\{ (x,y)\in \cB \ : \ y^2\equiv f(x) \pmod{q}\}|\le  \frac{H^{3/2}}{q^{1/6}}+H^{1/3+o(1)}.$$
In particular if $H\le q^{1/7}$ then
$$|\{ (x,y)\in \cB \ : \ y^2\equiv f(x) \pmod{q}\}|\le  H^{1/3+o(1)}.$$
\end{theorem}

For comparison with Theorem~\ref{thm:ConcentEllipticWiderRange}, we note that, for a prime $p$,  the bounds of~\cite[Theorem 1]{CCGHS} and~\cite[Theorem~2]{CCGHS} imply that

\begin{align}
\label{eqn:CCGHS-Cor3}
& \nonumber |\{ (x,y)\in \cB \ : \ y^2\equiv f(x) \pmod{p}\}|\le  \\  &  \quad \quad \quad \quad \quad \quad \quad \quad    \begin{cases}H^{1/3+o(1)}, \quad \text{if} \quad H<p^{1/8}, \\  \left(\frac{H^4}{p} \right)^{1/6}H^{1+o(1)}, \quad \text{if} \quad p^{1/8}\le H<p^{5/23}, \\ \left(\frac{H^3}{p}\right)^{1/16}H^{1+o(1)}, \quad \text{if} \quad p^{5/23}\le H < p^{1/3}.\end{cases}
\end{align}
We remark that one may incorporate Lemma \ref{lem:vmvt1}, with $(k, s) = (3, 6)$ instead of~\cite[Lemma 13]{CCGHS}, into the proof of~\cite[Theorem~2]{CCGHS}, to obtain the estimate
\begin{equation}
\label{eqn:12vs16}
 |\{ (x,y)\in \cB \ : \ y^2\equiv f(x) \pmod{p}\}|\le H^{1/3+o(1)}+\(\frac{H^3}{p}\)^{1/12} H^{1+o(1)}.
\end{equation}
This improves on the bound (\ref{eqn:CCGHS-Cor3}) in the range $p^{1/5}\le H < p^{1/3}$ and appears to be the best bound possible through the arguments of~\cite[Theorem~2]{CCGHS}.
We see that Theorem~\ref{thm:ConcentEllipticWiderRange} improves on the bounds (\ref{eqn:CCGHS-Cor3}) and (\ref{eqn:12vs16}) in the range $p^{1/8}\le H < p^{1/3}.$

\section*{Acknowledgement}
The authors would like to thank Igor Shparlinski for his comments and for pointing out an error in a previous version of the current paper.
\section{Preliminary results}
We first recall some facts from the geometry of numbers. Given a lattice $\Gamma \subset \R^{n}$ and a symmetric convex body $D\subset \R^{n}$ we define the $i$-th successive minimum of $\Gamma$ with respect to $D$ by 
$$\lambda_{i}=\inf\{ \lambda : \Gamma \cap \lambda D \ \ \text{contains $i$ linearly independent points} \}.$$
We define the dual lattice $\Gamma^{*}$ by
$$\Gamma^{*}=\{ y\in \R^{n} \ : \ \langle y,z\rangle\in \Z \ \ \text{for all $z\in \Gamma$}\},$$
and the dual body $D^{*}$ by
$$D^{*}=\{ y \in \R^{n} \ : \ \langle y,z \rangle \le 1 \ \ \text{for all $z\in D$} \},$$
where $\langle   .  ,  .   \rangle$ denotes the Euclidean inner product. The following is Minkowski's second theorem, for a proof see~\cite[Theorem~3.30]{TV}.
\begin{lemma}
\label{lem:mst}
Let $\Gamma\subset \R^{n}$ be a lattice, $D\subset \R^{n}$ a symmetric convex body and let $\lambda_{1},\dots,\lambda_{n}$ denote the successive minima of $\Gamma$ with respect to $D$. We have
$$\frac{\text{Vol}(D)}{\text{Vol}(\R^n/\Gamma)}\ll \frac{1}{\lambda_1\dots\lambda_n}  \ll\frac{\text{Vol}(D)}{\text{Vol}(\R^n/\Gamma)}.$$
\end{lemma}
 We may bound the number of lattice points $|\Gamma \cap D|$ in terms of the successive minima, see for example~\cite[Exercise~3.5.6]{TV}.
\begin{lemma}
\label{lem:lattice}
Let $\Gamma\subset \R^{n}$ be a lattice, $D\subset \R^{n}$ a symmetric convex body and let $\lambda_{1},\dots,\lambda_{n}$ denote the successive minima of $\Gamma$ with respect to $D$. We have
$$|\Gamma \cap D|\ll \prod_{j=1}^{n}\max\left(1, \frac{1}{\lambda_j} \right).$$
\end{lemma}
The successive minima of a lattice with respect to a convex body and the successive minima of the dual lattice with respect to the dual body are related through the following estimates, see for example~\cite{Ba}.
\begin{lemma}
\label{lem:transfer}
Let $\Gamma\subset \R^{n}$ be a lattice, $D\subset \R^{n}$ a symmetric convex body and let $\lambda_{1},\dots,\lambda_{n}$ denote the successive minima of $\Gamma$ with respect to $D$. Let $\Gamma^{*}$ denote the dual lattice, $D^{*}$ the dual body and let $\lambda_1^{*},\dots,\lambda_n^{*}$ denote the successive minima of $\Gamma^{*}$ with respect to $D^{*}$. For each $1\le j \le n$ we have
$$\lambda_j \lambda^*_{n-j+1}\ll 1.$$
\end{lemma} 
The following is a consequence of results of Bourgain, Demeter and Guth~\cite[Theorem~4.1]{BDG} and of Wooley~\cite[Theorem~1.1]{Wo}.
\begin{lemma}
\label{lem:vmvt}
Let $k$ be an integer, $s$ a positive real number with $s\le k(k+1)/2$ and $(a_n)_{n\in \Z}$ a sequence of complex numbers. We have 
\begin{align*}
\int_{[0,1)^{k}}\left|\sum_{|n|\le X}a_ne(\alpha_1 n+\dots+\alpha_k n^{k}) \right|^{2s}d\boldsymbol{\alpha} \ll X^{o(1)}\left(\sum_{|n|\le X}|a_n|^2 \right)^{s}.
\end{align*}
\end{lemma}
The following is an immediate corollary of Lemma~\ref{lem:vmvt}.
\begin{lemma}
\label{lem:vmvt1}
Let $\cX\subset [1,X]\cap \Z$ be some set. For integers $k$ and $s$ we let $J_{k,s}(\cX)$ denote the number of solutions to the system of equations
\begin{align*}
x_1^{j}+\dots+x_{s}^{j}=x_{s+1}^{j}+\dots+x_{2s}^{j}, \quad 1\le j \le k,
\end{align*}
with variables satisfying
$$x_1,\dots,x_{2s}\in \cX.$$
For $s\le k(k+1)/2$ we have 
$$J_{k,s}(\cX)\ll |\cX|^{s}X^{o(1)}.$$
\end{lemma}

\section{Proof of Theorem~\ref{thm:main1}}
Supposing the cube $\cB$ is given by
$$\cB=(K,K+H]\times (L,L+H],$$
 we may assume $K,L=0$ by applying shifts $x\rightarrow x-K$, $y\rightarrow y-L$ and modifying the coefficients of $f$. Hence writing
\begin{align}
\label{eq:fdef}
f(x)=a_0+a_1x+\dots+a_d x^{d},
\end{align}
and letting $N$ denote the number of solutions to the congruence
\begin{align}
\label{eq:polyeq}
y\equiv f(x) \pmod{q},
\end{align}
with 
\begin{align}
\label{eq:polycond}
1\le x \le H, \quad 1\le y \le H,
\end{align}
it is sufficient to show that 
$$N\le \frac{H^{1+2/d(d+1)+o(1)}}{q^{2/d(d+1)}}+H^{1/d+o(1)}.$$
We define the lattice 
\begin{align}
\label{eq:Ldef}
\cL=\{(z,w_1,\dots,w_d)\in \Z^{d+1} : \ z+a_1w_1+\dots+a_d w_d \equiv 0 \pmod{q} \},
\end{align}
 the convex body
\begin{equation}
\begin{split}
\label{eq:Ddef}
D= &\{(h_0,\dots,h_d)\in \R^{d+1} \ : \\ &|h_0|\le dH \ \ \text{and} \ \  |h_i|\le dH^{i} \ \ \text{for} \ \  1\le i \le d \},
\end{split}
\end{equation}
and let $\lambda_1,\dots,\lambda_{d+1}$ denote the successive minima of $\cL$ with respect to $D$. We distinguish two cases depending on whether $\lambda_{d+1}<1$ or not.

 We first consider when 
\begin{align}
\label{eq:case1}
\lambda_{d+1}<1.
\end{align}
 Let 
$$s=\frac{d(d+1)}{2},$$
and suppose
$(x_1,y_1),\dots,(x_{2s},y_{2s})$ are $2s$ solutions to~\eqref{eq:polyeq} satisfying~\eqref{eq:polycond}.
We have 
\begin{align}
\label{eq:addsols}
z\equiv f(x_1)+\dots+f(x_{s})-f(x_{s+1})-\dots-f(x_{2s}) \pmod{q},
\end{align}
for some $|z|\le sH.$ We define the set $\cX$ by
\begin{align*}
&\cX=\\ &\{1\le x \le H : \  \text{there exists} \ \  1\le y \le H \ \ \text{such that} \ \ y\equiv f(x) \Mod{q}\},
\end{align*}
so that
\begin{align}
\label{eq:Xcard}
|\cX|=N.
\end{align}
Let $J(w_1,\dots,w_{d})$ denote the number of solutions to the system of equations
\begin{align}
x_1^{j}+\dots+x_{s}^{j}-x_{d+1}^{j}-\dots-x_{2s}^{j}\equiv w_j \pmod{q}, \quad 1\le j \le d,
\end{align}
with variables  satisfying 
$$x_1,\dots,x_{2s}\in \cX,$$
and when $$w_1=\dots=w_{d}=0,$$ we write
$$J(0,\dots,0)=J.$$
The equation~\eqref{eq:addsols} implies that 
\begin{align*}
N^{2s}\le \sum_{|z|\le sH}\sum_{\substack{|w_i|\le sH^{i} \\ z\equiv a_1w_1+\dots+a_dw_d \pmod{q} }}J(w_1,\dots,w_d),
\end{align*}
and since
\begin{align*}
J(w_1,&\dots,w_d) \\ &=\int_{(0,1]^{d}}\left|\sum_{x\in \cX}e^{2\pi i(\alpha_1 x+\dots+\alpha_d x^{d})} \right|^{2d}e^{-2\pi i(\alpha_1w_1+\dots+\alpha_d w_d)}d\alpha_1\dots d\alpha_{d} \\ 
&\le \int_{(0,1]^{d}}\left|\sum_{x\in \cX}e^{2\pi i(\alpha_1 x+\dots+\alpha_d x^{d})} \right|^{2d}d\alpha_1\dots d\alpha_{d}=J,
\end{align*}
the above gives
\begin{align}
\label{eq:N2din}
N^{2s}\le \left(\sum_{|z|\le dH}\sum_{\substack{|w_i|\le dH^{i} \\ z\equiv a_1w_1+\dots+a_dw_d \pmod{q} }}\right)J\le |\cL \cap D|J,
\end{align}
where $\cL$ and $D$ are given by~\eqref{eq:Ldef} and~\eqref{eq:Ddef}.
We recall that $J$ denotes the number of solutions to the system
$$x_1^{j}+\dots+x_{s}^{j}=x_{d+1}^{j}+\dots+x_{2s}^{j}, \quad 1\le j \le d,$$
with 
$$x_1,\dots,x_{2s}\in \cX.$$
By Lemma~\ref{lem:vmvt1} we have 
$$J\ll |\cX|^sH^{o(1)}= H^{o(1)}N^{s},$$
which by~\eqref{eq:N2din} implies that
\begin{align}
\label{eq:Nlattice}
N^{s}\ll |\cL\cap D|H^{o(1)}.
\end{align}
Combining the assumption~\eqref{eq:case1} with Lemma~\ref{lem:mst} and Lemma~\ref{lem:lattice} gives
\begin{align*}
|\cL \cap D|\ll \frac{1}{\lambda_1\dots \lambda_{d+1}}\ll \frac{\text{Vol}(D)}{\text{Vol}(\R^{d+1}/\cL)}.
\end{align*}
Since 
\begin{align*}
\text{Vol}(D)\ll H^{(d^2+d+2)/2} \quad \text{and} \quad \text{Vol}(\R^{d+1}/\cL)=q,
\end{align*}
we see that 
$$|\cL\cap D|\ll \frac{H^{(d^2+d+2)/2}}{q},$$
and hence by~\eqref{eq:Nlattice} 
\begin{align}
\label{eq:case1final}
N\ll H^{o(1)}\left(\frac{H^{(d^2+d+2)/2}}{q}\right)^{1/s}=\frac{H^{1+2/d(d+1)+o(1)}}{q^{2/d(d+1)}}.
\end{align}

Consider next when 
\begin{align}
\label{eq:case2}
\lambda_{d+1}\ge 1.
\end{align}
The dual lattice $\cL^{*}$  and the dual body $D^{*}$ are given by
\begin{align*}
\cL^{*}= \frac{1}{q}&\{(z,w_1,\dots,w_d) \in \Z^{5}~:\\  &\exists \  n\in \Z \ \ \text{such that} \ \   n\equiv z \pmod{q},\ \ a_i n\equiv w_i \pmod{q}\},
\end{align*}  
and
\begin{align*}
D^{*}=\{ (h_0,h_1,\dots,h_d)\in \R^{5}~: sH|h_0|+\sum_{i=1}^{d}sH^{i}|h_i|\le 1 \}. 
\end{align*}

If $\lambda_1^{*}$ denotes the first successive minimum of $\cL^{*}$ with respect to $D^{*}$ then by~\eqref{eq:case2} and Lemma~\ref{lem:transfer} there exists some constant $c$ depending only on $d$ such that
\begin{align*}
\lambda_1^{*}\le c.
\end{align*}
This implies that 
$$\cL^{*}\cap cD^{*}\neq \emptyset,$$
and hence there exist $z,w_1,\dots,w_d,n\in \Z$ such that
\begin{align}
\label{eq:congcond}
n\equiv z \pmod{q}, \quad a_in\equiv w_i \pmod{q},
\end{align}
and
\begin{align}
\label{eq:sizecond}
|z|\ll \frac{q}{H}, \quad |w_i|\ll \frac{q}{H^{i}}.
\end{align}
We note that since $(a_d,q)=1$ we have $w_d\neq 0$. Supposing $(x,y)$ is a solution to~\eqref{eq:polyeq} and recalling~\eqref{eq:fdef}, we have
$$ny\equiv a_0n+a_1nx+\dots+a_dnx^{d} \pmod{q},$$
so that writing $w_0$ for the least residue of $a_0n \pmod{q}$, an application of~\eqref{eq:congcond} gives
\begin{align*}
w_0+w_1x+\dots+w_dx^{d}-zy=tq,
\end{align*}
for some $t\in \Z$. By~\eqref{eq:sizecond}, the number of possible values for $t$ is $O(1)$ and for each such value of $t$ we apply Theorem~\ref{thm:bombieripila}. This gives
\begin{align}
\label{eq:case2final}
N\le H^{1/d+o(1)},
\end{align}
and the result follows by combining~\eqref{eq:case1final} and~\eqref{eq:case2final}.
\section{Proof of Theorem~\ref{thm:ConcentEllipticWiderRange}}
Applying shifts as in the proof of Theorem~\ref{thm:main1}, it suffices to estimate the number of solutions to a congruence of the form
\begin{equation}
\label{eq:ellipticAfterShifts}
y^2-c_0y\equiv a_3x^3+a_2x^2+a_1x+a_0\pmod q,\quad 1\le x,y\le H,
\end{equation}
where $(a_3,q)=1$.
Let $N$ denote the number of solutions to the congruence~\eqref{eq:ellipticAfterShifts} and let $\cX$ denote the set of $x$ for which $(x,y)$ satisfies~\eqref{eq:ellipticAfterShifts} for some $1\le y \le H$, so that
\begin{equation}
\label{eq:NX}
N\le  2|\cX|.
\end{equation}
For $j=1,2,3$ we define the intervals
$$
I_{j}=[-6 H^j, 6 H^j],
\qquad  j =1,2,3,
$$
and consider the set
$$
\cS\subseteq  I_{1}\times I_{2}\times I_{3},
$$
of all triples
\begin{equation}
\label{eq:triples} \vec{x}=
(x_1+\ldots +x_{6}, \, x_1^2+\ldots
+x_{6}^2, \, x_1^3+\ldots +x_{6}^3),
\end{equation}
such that $x_i\in \cX$. For each $\vec{x}$ we let $I(\vec{x})$ count the number of solutions to the equation~\eqref{eq:triples} with variables $x_1,\dots,x_6\in \cX$. Thus
$$\sum_{\vec{x} \in \cS} I(\vec{x}) = |\cX|^{6},$$
and 
$$\sum_{\vec{x} \in \cS} I(\vec{x})^2,$$
 is bounded by the number of solutions
to the system of equations
\begin{equation}
\label{eq:congrWooley} x_1^j+\ldots +x_{6}^j\equiv x_{7}^j+\ldots
+x_{12}^j,\qquad j=1,2,3,
\end{equation}
with variables 
$$x_1,\dots,x_{12}\in \cX.$$
By Lemma~\ref{lem:vmvt1}
$$\sum_{\vec{x} \in \cS} I(\vec{x})^2\le H^{o(1)}|\cX|^{6},$$
and by the Cauchy-Schwarz inequality, we have 
$$
\sum_{\vec{x} \in \cS} I(\vec{x}) \le \(|\cS|
\sum_{\vec{x} \in \cS} I(\vec{x})^2\)^{1/2},
$$
which implies that
$$
|\cS|\ge |\cX|^{6}H^{o(1)}.
$$
Hence, there exist
at least $|\cX|^{6}H^{o(1)}$ triples
$$
(x_1,x_2,x_3)\in I_{1} \times I_{2} \times I_{3},
$$
such that
$$
a_3x_3+a_2x_2+a_1x_1\equiv y_2-c_0y_1\pmod q,
$$
for some $y_2\in I_{2}$ and $y_1\in I_{1}$.
In particular, we have that the congruence
\begin{equation*}
\begin{split}
a_1x_1+a_2x_2+a_3x_3 & +c_0y_1+y_2\equiv 0\pmod q,\\
(x_1,x_2,x_3,y_1,y_2)&\in I_{1}\times I_{2} \times I_{3}\times I_{1}\times I_{2},
\end{split}
\end{equation*}
has a set of solutions $\cS$ with
\begin{equation}
\label{eq:S large}
| \cS|\ge |\cX|^{6}H^{o(1)}.
\end{equation}


We define the lattice
\begin{equation*}
\begin{split}
\cL = \{(x_1,x_2,x_3,y_1,y_2)&\in\Z^5~:\\
~
a_1x_1+a_2x_2+a_3x_3 & +c_0y_1+y_2\equiv 0\pmod q\},
\end{split}
\end{equation*}
and the body
\begin{equation*}
\begin{split}
D = &\{(x_1,x_2,x_3,y_1,y_2)\in\R^5~:\\ &|x_1|, |y_1|\le 6 H, \ |x_2|, |y_2|\le 6 H^{2},\ |x_3|\le 6 H^{3}\}.
\end{split}
\end{equation*}
It follows from~\eqref{eq:S large} that
\begin{align}
\label{eq:Xl1}
|\cX|^{6}\ll |\Gamma \cap D|H^{o(1)}.
\end{align}
Let $\lambda_i$ denote the $i$-th successive minimum of  $\cL$ with respect to $D$.
 We consider two cases depending on whether $\lambda_5\le 1$ or not. Suppose first that $\lambda_5\le 1$.  By Lemma~\ref{lem:mst} and Lemma~\ref{lem:lattice}
\begin{align*}
|\Gamma \cap D|\ll \frac{H^9}{q},
\end{align*}
and hence by~\eqref{eq:Xl1}
\begin{align}
\label{eq:Xc1}
|\cX|\le \frac{H^{3/2+o(1)}}{q^{1/6}},
\end{align}
and the result follows from~\eqref{eq:NX}.
Suppose next that
\begin{align}
\label{eq:lambda5in}
\lambda_5>1.
\end{align}
The dual lattice $\cL^{*}$ and dual body $D^{*}$ are given by 
\begin{equation*}
\begin{split}
 \cL^{*}=\frac{1}{q}&\{ (w_1,w_2,w_3,z_1,z_2)\in \Z^{5}~:\\ &\exists \ n\in \Z \ \ \text{such that} \ \  w_i\equiv a_in  \pmod{q},\  i=1,2,3,\\ &z_1\equiv c_0 n \pmod{q}, \ \ \text{and} \ \  \ z_2\equiv n \pmod{q} \}.
\end{split}
\end{equation*}
and
\begin{align*}
D^{*}=\{ (w_1,w_2,w_3,z_1,z_2) \ : \ \sum_{i=1}^{3}6H^{i}|w_i|+6H|z_1|+6H^2|z_2|\le 1 \}.
\end{align*}
If $\lambda_1^{*}$ denotes the first successive minimum of $\cL^{*}$ with respect to $D^{*}$, by Lemma~\ref{lem:transfer} we have 
$$\lambda_5 \lambda_1^{*} \ll 1,$$
and hence by~\eqref{eq:lambda5in}
$$\lambda_1^{*}\ll 1.$$
This implies that there exists some integer $n$ satisfying
\begin{align}
\label{eq:nbound1}
n\ll \frac{q}{H^2},
\end{align}
and
\begin{align*}
a_in \equiv w_i \pmod{q}, \ \  i=1,2,3 \quad \text{and} \quad  c_0n \equiv z_1 \pmod{q},
\end{align*}
for some $w_1,w_2,w_3$ and $z_1$ satisfying
\begin{align}
\label{eq:wz}
w_i \ll \frac{q}{H^i}, \ \ \ i=1,2,3, \quad z_1 \ll \frac{q}{H}.
\end{align}
For any solution $(x,y)$ to~\eqref{eq:ellipticAfterShifts} we  have
$$ny^2-z_1y=w_3x^3+w_2x^2+w_1 x+qt,$$
for some $t\in \Z$. By~\eqref{eq:nbound1} and~\eqref{eq:wz} there are  $O(1)$ possible values of $t$ and for each such value of $t$ we apply Theorem~\ref{thm:bombieripila}. This implies that the number of solutions to the congruence~\eqref{eq:ellipticAfterShifts} is bounded by
 $$N\le H^{1/3+o(1)},$$
and the result follows combining the above with~\eqref{eq:NX} and~\eqref{eq:Xc1}.

\end{document}